\documentclass[11pt]{amsart}
\usepackage{amscd,amsmath,amssymb,amsfonts,verbatim}
\usepackage[cmtip, all]{xy}
\usepackage{comment}
\usepackage{tikz-cd}
\usepackage{hyperref}

\renewcommand\qedsymbol{$\blacksquare$}

\newtheorem{thm}{Theorem}[section]

\newtheorem{cor}[thm]{Corollary}

\theoremstyle{definition}

\newtheorem{defn}[thm]{Definition}
\newtheorem{remks}[thm]{Remarks}
\newtheorem{remk}[thm]{Remark}

\newtheorem{exm}[thm]{Example}
\newtheorem{exms}[thm]{Examples}
\newtheorem{notat}[thm]{Notation}
\numberwithin{equation}{section}
\numberwithin{equation}{section}

{\hfill$\square$\end{defn}}
{\end{remk}}
{\hfill$\square$\end{remks}}
{\hfill$\square$\end{exm}}
{\hfill$\square$\end{exms}}
{\hfill$\square$\end{notat}}

\newcommand{\A}{\mathrm{A}}
\newcommand{\C}{\mathrm{C}}
\newcommand{\D}{\mathrm{D}}
\newcommand{\E}{\mathrm{E}}
\newcommand{\G}{\mathrm{G}}

\newcommand{\rank}{{\rm rank}}

\newcommand{\Max}{{\rm Max \ }}

\newcommand{\Ker}{{\rm Ker}}
\newcommand{\Img}{{\rm Im}}

\newcommand{\El}{{\operatorname{\rm E}}}

\newcommand{\K}{{\operatorname{\rm K}}}

\newcommand{\St}{\mathrm{St}}
\newcommand{\ev}{\mathrm{ev}}

\newcommand{\ds}{{/\kern-3pt/}}

\newcommand{\tuborg}{\left\{\begin{array}{ll}}
\newcommand{\sluttuborg}{\end{array}\right.}

\newcounter{elno}

\newcounter{elno-abc}   

\newcounter{elno-abc-prime}

\begin{document}

\title[Graded version of local global principle]{Graded version of local-global principle for Steinberg and Chevalley groups}
\author{Chandan Bhaumik}

\address{Department of Mathematics, Indian  Institute of Technology Bombay, Powai, Mumbai 400076, India}
 \email{cbhaumik11math@gmail.com, cbhaumik@math.iitb.ac.in}

\dedicatory{}

\keywords{Steinberg group, Chevalley group, Local–global principle, Graded ring}

\subjclass[2020]{20G35, 19C20, 13A02}

\maketitle
\vspace{-.7cm}

\begin{quote}\emph{Abstract.} 
In this article, an analog of a graded version of Quillen–Suslin’s local-global principle was established for Steinberg groups $\St(\Phi,-)$ and Chevalley groups $\G(\Phi,-)$ over a reduced irreducible root system $\Phi$.  
\end{quote}
\setcounter{tocdepth}{1} 


\section{Introduction}\label{sec:Intro}
In the context of classical K-theory, the local-global principle (LGP for short) was introduced by Quillen \cite{Qu76} and Suslin \cite{Su77} as part of their solution to Serre’s problem on projective modules over polynomial extensions of a field. This principle was initially applied to prove Serre’s conjecture. The solution to Serre’s problem was accompanied by the development of several techniques that are interesting in their own right, including various versions of LGP.\par
Recall that for a commutative ring $A$, denote $\G(\Phi,A)$ as a Chevalley-Demazure group scheme (see \cite{St68}) and $\St(\Phi,A)$ as a Steinberg group presented by generators $x_{\alpha}(a)$ for $\alpha\in\Phi$ and $a\in A$ (see Definition $\ref{def_st_grp}$) with a reduced irreducible root system $\Phi$. The subgroup of $\G(\Phi,A)$ generated by the elementary root unipotents $e_{\alpha}(a)$ is called the elementary subgroup of $\G(\Phi,A)$ and is denoted by $\El(\Phi,A)$. There is a well-defined map $\phi:\St(\Phi,A)\to\G(\Phi,A)$ sending $x_{\alpha}(a)\mapsto e_{\alpha}(a)$, and the image of $\phi$ is $\El(\Phi,A)$, which is a normal subgroup of $\G(\Phi,A)$ when $\rank(\Phi)\geq 2$. The groups $\K_{1}(\Phi,A)$ and $\K_{2}(\Phi,A)$ are defined as the cokernel and kernel of $\phi$, respectively
\begin{center}
        \begin{tikzcd}
1 \arrow[r] & \K_{2}(\Phi,A) \arrow[r]           & \St(\Phi, A) \arrow[r,"\phi"]                 & \G(\Phi, A) \arrow[r]  & \K_{1}(\Phi,A) \arrow[r]               & 1. 
\end{tikzcd}
\end{center}\par
The $\K_{1}$-analog of the LGP was established for Chevalley groups of type $\A_{n}$ by Suslin \cite[Theorem 3.1]{Su77}. Subsequent work extended it to type $\D_{n}$ by Suslin–Kopeiko \cite[Theorem 4.2]{SuKo82}, to type $\C_{n}$ by Grunewald–Mennicke–Vaserstein \cite[Proposition 5.6]{GMV91}, and to all  Chevalley groups by Abe \cite[Theorem 1.15]{Ab83} and Wendt \cite[Corollary 4.3]{We10}. Turning to the $\K_{2}$-analog of the LGP, it was proven for Chevalley groups of type $\A_{n}$ by Tulenbaev \cite[Theorem 2.1]{Tu83}. Recently, it was extended to type $\E_{n}$ by Sinchuk \cite[Theorem 2]{Si16}, to type $\D_{n}$ by Lavrenov-Sinchuk \cite[Theorem 2]{LaSi17}, and to type $\C_{n}$ by Lavrenov \cite{La18}.\par
The centrality of $\K_{2}$ (see \cite{LaSi17}) follows from the $\K_{2}$-analog of Quillen’s LGP, just as the normality theorem for the elementary subgroup is derived from the $\K_{1}$-analog of Suslin's LGP. In other words, the LGP plays a foundational role in establishing both the normality of $\El(\Phi,A)$ and the centrality of $\St(\Phi,A)$. The present article aims to revisit these results in the setting of commutative nonnegative graded rings.\par
The graded version of the $\K_{1}$-analog of the LGP has been studied for classical groups in \cite{BaSi20,Gu92}. Meanwhile, for Chevalley groups and Steinberg groups, the graded versions of both the $\K_{1}$- and $\K_{2}$-analogs of the LGP are expected to hold. We establish the general statement (Theorem $\ref{thm_lgptoglgp}$), namely the graded versions of the LGP (GLGP for short) for a group-valued functor, under the assumption that the corresponding LGP is already known. The main technical ingredient in the proof is the Swan--Weibel homotopy trick, which enables us to pass from the result for polynomial rings to the corresponding result for nonnegatively graded rings.\par
As a consequence of the main theorem ($\ref{thm_lgptoglgp}$), we therefore establish the following GLGP for Steinberg groups and Chevalley groups, respectively.
\begin{thm}\label{intro_thm_lgp_stgp}
    Let $R=\bigoplus_{i\geq 0}R_{i}$ be a commutative graded ring and $\Phi_{n}$ be either a root system of type simply-laced or $\C_{n}$ of rank $\geq 3$. Suppose $\theta\in \St(\Phi_{n},R, R^{+})$. Then \begin{enumerate}
        \item $\theta=1$ if and only if $\lambda^{\ast}_{\mathfrak{m}}(\theta)=1\in \St(\Phi_{n},R_{\mathfrak{m}})$ for all maximal ideals $\mathfrak{m}$ of $R_{0}$.
        \item $\theta\in\Img(l_{R})$ if and only if $\lambda^{\ast}_{\mathfrak{m}}(\theta)\in \Img(l_{R_{\mathfrak{m}}})$ for all maximal ideals $\mathfrak{m}$ of $R_{0}$.
    \end{enumerate}
\end{thm}
\begin{thm}\label{intro_thm_lgp_chgp}
     Let $R=\bigoplus_{i\geq 0}R_{i}$ be a commutative graded ring and $\Phi$ be a root system of rank $\geq 2$. Assume that $R$ has no residue fields of $2$ elements if $\Phi=\C_{2}~\text{or}~\G_{2}$. Suppose $\theta\in \G(\Phi,R,R^{+})$. Then $\theta\in \El(\Phi,R)$ if and only if $\lambda^{\ast}_{\mathfrak{m}}(\theta)\in \El(\Phi,R_{\mathfrak{m}})$ for all maximal ideals $\mathfrak{m}$ of $R_{0}$.
\end{thm}

\section{Preliminaries and Notations}\label{sec_preli}
Let $A$ be a commutative ring with $1$ and $\Phi$ being a reduced irreducible root system.
\begin{defn}\label{def_st_grp}
The {\it Steinberg group} $\St(\Phi, A)$ is the group generated by symbols $x_\alpha(u)$, where $\alpha \in \Phi$ and $u \in A$, subject to the following relations:
\begin{align*}
    &    x_\alpha(u) \, x_\alpha(v) = x_\alpha(u+v), \qquad u, v \in A,~ \alpha\in\Phi;\\[4pt]
    &    [x_\alpha(u), x_\beta(v)] = \prod_{\substack{i\alpha + j\beta \in \Phi}} x_{i\alpha + j\beta}\bigl(C_{i,j,\alpha,\beta}\, u^i v^j\bigr), \qquad u, v \in A,~\alpha, \beta \in \Phi, ~\alpha\neq -\beta,
\end{align*}
where $[p,q] = pqp^{-1}q^{-1}$, the product is taken over all positive integers $i,j$ such that $i\alpha + j\beta \in \Phi$, and $C_{i,j,\alpha,\beta} \in \mathbb{Z}$ are the Chevalley structure constants (depending only on $\Phi$) of the simply connected Chevalley group $\G(\Phi,A)$. More information on Chevalley groups, the structure constants, and various types of root systems $\Phi$ can be found in \cite{St68,VaPl96}. 
\end{defn}
For an ideal $I$ of $A$, the subgroup $\St(\Phi,I)$ is spanned by symbols $x_\alpha(s)$, where $\alpha \in \Phi$ and $s \in I$. An ideal $I$ is called a {\it splitting ideal} if the canonical projection ring homomorphism $\pi:A \twoheadrightarrow A/I$ splits, i.e., there exists $\sigma:A/I\to A$ such that $\pi\sigma=1$. It is noted that $TA[T]$ is a splitting ideal of $A[T]$.
\begin{defn}
  For a splitting ideal $I$, the {\it relative Steinberg group} $\St(\Phi,A,I)$ is defined as the kernel  $\Ker(\pi^{\ast}:\St(\Phi,A)\twoheadrightarrow\St(\Phi,A/I))$. Obviously, $\St(\Phi,A,I)$ is coincides with the normal closure of $\St(\Phi,I)$ in $\St(\Phi,A)$. The congruence subgroup $\G(\Phi,A,I)$ is the kernel of the canonical map $\G(\Phi,A)\to \G(\Phi,A/I)$.  
\end{defn}
   
\begin{defn}
    A root system $\Phi$ is called {\it simply-laced} if the Dynkin diagram of $\Phi$ does not contain double bonds (i.e., $\Phi$ has type $\A_n,\D_n$ and $\E_{6,7,8}$). 
\end{defn}
\textbf{Notations.} For a ring $A$, $A_{\mathfrak{m}}$ denotes the localization at a maximal ideal $\mathfrak{m}$ of $A$. For a graded ring $R=\bigoplus_{i\geq 0}R_{i}$, $R_{\mathfrak{m}}$ denotes the localization at a maximal ideal of its degree-zero component $R_0$. The localization homomorphism $A\to S^{-1}A$ is denoted by $\lambda_{S}$ for a multiplicative closed subset $S$ of $A$. If $f:A\to B$ is a ring homomorphism and $G$ is a covariant functor from the category of rings to the category of groups, its corresponding group homomorphism is denoted by $f^{\ast}:\G(A)\to \G(B)$. If $\Phi_{n}$ is an irreducible root system of rank $n$ and $\Phi_{n-1}$ is the subsystem of same type, we denote by $l_{A}:\St(\Phi_{n-1},A)\to \St(\Phi_{n},A)$ the homomorphism of Steinberg groups induced by a map of root systems $\Phi_{n-1}\to \Phi_n$. 
\section{Graded Version of Local-Global Principle}\label{sec:lgp}
One of the key components for handling the graded situation is the well-known ``Swan--Weibel homotopy trick".
Let $R=R_{0}\oplus R_{1}\oplus R_{2}\oplus\cdots$ be a commutative graded ring and $R^{+}=R_{1}\oplus R_{2}\oplus\cdots$ denote the irrelevant ideal of $R$. It is noted that $R^+$ is the splitting ideal of $R$. An element $r\in R$ will be denoted by $r=r_0+r_1+r_2+\cdots$ where $r_i\in R_i$ for each $i$, and all but finitely many $r_{i}$ are zero.
\begin{defn}\label{def_sw_map}
    For a graded ring $R=\bigoplus_{i\geq 0}R_{i}$, the {\it Swan--Weibel homotopy map} $h:R\to R[T]$ is defined as \[h(r_0+r_1+\cdots+r_m)=r_0+r_1 T+\cdots r_m T^{m}, ~(r_{i}\in R_{i}).\]
For any $R$-algebra $A$ and any $a\in A$, consider an evaluation map $\ev_{a}:R[T]\to A$, which is the $R$-algebra homomorphism sending $T\mapsto a$. For more details on Swan--Weibel homotopy trick, we refer the reader to \cite[\S V.3]{Lambook}.
\end{defn}
We denote the image $r$ in $A$ by $r^{h}(a):=(\ev_{a}\circ h)(r)$. If $A=R$, then it is easy to observe that $\ev_{0}\circ h: R\to R$ is precisely the projection onto the degree-zero part and $\ev_{1}\circ h: R\to R$ is precisely the identity map. For a fixed element $c\in R_0$, similarly denote $r^{h}(c):=h(r)(c)= (\ev_{c}\circ h)(r)$, and note that $r^{h}(0)=r_0$ and $r^{h}(1)=r$. In order to emphasize the choice of $A$, we denote by $\ev\left [\frac{R[T]\to A}{T\mapsto a}\right ]$ for $R$-algebra homomorphism $\ev_{a}$. \par
Let $\G$ be a group-valued covariant functor and $\El$ be a group-valued subfunctor from the category of commutative rings,
\[\El\subset\G:\text{CommRings}\to\text{Groups}.\]
The map $h:R\to R[T]$ (resp. $\ev_{a}:R[T]\to A$) induces a group homomorphism $h^{\ast}:\G(R)\to \G(R[T])$ (resp. $\ev_{a}^{\ast}:\G(R[T])\to \G(A)$) since $\G$ is a covariant functor. For $\theta\in \G(R)$ and $c\in R_{0}$, by functorial properties, analogously, $\theta^{h^{\ast}}(c):=(\ev_{c}^{\ast}\circ h^{\ast})(\theta)$ is defined. \par
\textbf{Local-Global Principle.} The ordinary LGP for $(\El,\G)$ states that, for an arbitrary commutative ring $A$ and $\theta(X)\in\G(A[X])$ with $\theta(0)=1$, the conditions \[\lambda^{*}_{M}(\theta(X))\in\El(A_{M}[X]) \quad \forall~ M\in\max(A)\quad\text{implies} \quad\theta(X)\in \El(A[X]).\]\par
\smallskip
The following result is a graded version of the local-global principle for $(\El,\G)$ deduced from the LGP by Swan–Weibel homotopy trick. 
\begin{thm}\label{thm_lgptoglgp}
    Let $R=\bigoplus_{i\geq 0}R_{i}$ be a graded ring and `overbar' denote the reduction modulo the ideal $R^{+}$. Suppose $\theta\in\G(R)$ with $\overline{\theta}=1\in\G(R_0)$, Assume that $(\El,\G)$ satisfies LGP. Then \[\lambda^{*}_{\mathfrak{m}}(\theta)\in\El(R_{\mathfrak{m}}) \quad \forall~ \mathfrak{m}\in\Max(R_0)\quad\text{implies} \quad\theta\in \El(R).\]
\end{thm}
\begin{proof}
    Denote $\Gamma(X)=\theta^{h^{\ast}}(X)=\left(\ev\left [\frac{R[T]\to R[X]}{T\mapsto X}\right ]^{\ast}\circ h^{\ast}\right)(\theta)\in \G(R[X])$. For $M\in\Max(R)$, let $\mathfrak{p}=M\cap R_{0}$ and choose $\mathfrak{m}\in\Max(R_0)$ such that $\mathfrak{p}\subset \mathfrak{m}$. Since $\lambda^{*}_{\mathfrak{m}}(\theta)\in\El(R_{\mathfrak{m}})$, we have \[\lambda^{\ast}_{\mathfrak{m}}(\Gamma(X))=\left(\ev\left [\frac{R_{\mathfrak{m}}[T]\to R_{\mathfrak{m}}[X]}{T\mapsto X}\right ]^{\ast}\circ h_{\mathfrak{m}}^{\ast}\right)(\lambda^{\ast}_{\mathfrak{m}}(\theta))\in\El(R_{\mathfrak{m}}[X]).\] Since $R_{M}$ is obtained by further localizing $R_{\mathfrak{m}}$, we conclude that $\lambda^{\ast}_{M}(\Gamma(X))\in\El(R_{M}[X])$. Note that $\Gamma(0)=\theta^{h^{\ast}}(0)=\overline{\theta}=1$. The LGP for $(\El,\G)$ implies that $\Gamma(X)\in\El(R[X])$. Hence, substituting $X=1$ yields $\theta=\Gamma(1)\in \El(R)$.
\end{proof}
\begin{cor}
     Let $R=\bigoplus_{i\geq 0}R_{i}$ be a graded ring and $\theta\in\G(R)$. Assume that $(\El,\G)$ satisfies LGP. Then \[\lambda^{*}_{\mathfrak{m}}(\theta)\in\G((R_0)_{\mathfrak{m}})\cdot\El(R_{\mathfrak{m}}) \quad \forall~ \mathfrak{m}\in\max(R_0)\quad\text{implies} \quad\theta\in \G(R_0)\cdot\El(R).\]
\end{cor}
\begin{proof}
   Let `overbar' denote the reduction modulo the ideal $R^{+}$. After left-multiplying $\theta$ by $\overline{\theta}^{-1}$, we can assume that $\theta^{h^{\ast}}(0)=1$. For all maximal ideals $\mathfrak{m}$ of $R_{0}$, we have the decomposition $\lambda^{\ast}_{\mathfrak{m}}(\theta)=\eta\cdot \gamma$ where $\eta\in \G((R_{0})_{\mathfrak{m}})$ and $\gamma\in \El(R_{\mathfrak{m}})$. Now, reduction modulo $R^{+}$, we get $1=\eta\cdot\overline{\gamma}$, and hence $\eta=\overline{\gamma}^{-1}$. Thus $$\lambda^{\ast}_{\mathfrak{m}}(\theta)=\overline{\gamma}^{-1}\cdot\gamma\in\El(R_{\mathfrak{m}}).$$ Then, by Theorem $\ref{thm_lgptoglgp}$, we get $\theta\in\El(R)$. Therefore $\theta\in \G(R_{0})\cdot\El(R)$ since the original $\theta$ has already been modified by a left factor from $\G(R_{0})$.
\end{proof}
\textbf{Proof of Theorem $\ref{intro_thm_lgp_stgp}$.} $(1)$ Choose a group functor $\G=\St(\Phi_n,R)$ and $\El=1$. By \cite[Theorem 2]{LaSi17} for $\Phi_n=\A_{n\geq 3},\D_{n\geq 3},\E_{n=6,7,8}$ and \cite[Main Theorem]{La18} for $\Phi_n=\C_{n\geq 3}$, the pair $(\El,\G)$ satisfies LGP. \par
$(2)$ Choose a group functor $\G=\St(\Phi_n,R)$ and $\El=\Img(l_R)$. By \cite[Corollary 5.3]{Si25} for $\Phi_n=\A_{n\geq 3},\D_{n\geq 3},\E_{n=6,7,8}$ and \cite[Main Theorem]{La18} for $\Phi_n=\C_{\geq 3}$, the pair $(\El,\G)$ satisfies LGP.\par
Therefore, both $(1)$ and $(2)$, follow from Theorem $\ref{thm_lgptoglgp}$.\hfill\qedsymbol
\par\medskip
\textbf{Proof of Theorem $\ref{intro_thm_lgp_chgp}$.} Choose a group functor $\G=\G(\Phi,R)$ and $\El=\El(\Phi,R)$. By \cite[Theorem 7.1]{ApSt13} for a given root system $\Phi$, the pair $(\El,\G)$ satisfies LGP, and hence the result follow from Theorem $\ref{thm_lgptoglgp}$.\hfill\qedsymbol
\par\medskip
\textbf{Acknowledgments:} The author sincerely thanks the referee for detailed comments and suggestions for improvement in the exposition of the paper. The author thanks Prof. Manoj Kumar Keshari for suggesting the problem. This material is based on work supported by the IIT Bombay Postdoctoral Research Fellowship.


\begin{thebibliography}{99}
\bibitem{Ab83}
Abe, E.: Whitehead groups of Chevalley groups over
polynomial rings. Comm. Algebra {\bf 11}(12), 1271--1307 (1983)
\bibitem{ApSt13}
Apte, H., Stepanov, A.: Local-global principle for congruence subgroups of Chevalley groups. Central Eur. J. Math. {\bf 12}, 801--812 (2014)

\bibitem{BaSi20}
Basu, R., Singh, M.K.: Quillen-Suslin theory for classical groups: revisited over graded rings. Contemp. Math. {\bf 751}, 5--18 (2020)
\bibitem{GMV91} 
Grunewald, F., Mennicke, J., Vaserstein, L.N.: On symplectic groups over polynomial rings. Math. Zeitschrift {\bf 206}, 35--56 (1991)

\bibitem{Gu92}
Gubeladze, J.: The elementary action on unimodular rows over a monoid ring. J. Algebra {\bf 148}, 135--161 (1992)

\bibitem{Lambook}
Lam, T.Y.: Serre's problem on projective modules. Springer Monographs in Mathematics, Springer-Verlag, Berlin (2006)
\bibitem{La18}
Lavrenov, A.: A local–global principle for symplectic $\K_{2}$. Doc. Math. {\bf 23}, 653--675 (2018)
\bibitem{LaSi17} 
Lavrenov, A., Sinchuk, S.: On centrality of even orthogonal $\K_{2}$. J. of Pure and Appl. Algebra {\bf 221}(2), 1134--1145 (2017)

\bibitem{Qu76}
Quillen, D.: Projective modules over polynomial rings. Invent. Math. {\bf 36}, 167--171 (1976)

\bibitem{Si16}
Sinchuk, S.: On centrality of $\K_{2}$ for Chevalley groups of type $\E_{l}$. J. of Pure and Appl. Algebra {\bf 220}(2), 857--875 (2016)
\bibitem{Si25}
Sinchuk, S.: On the $\mathbb{A}^{1}$-invariance of $\K_{2}$ of Chevalley groups of simply-laced type. Eur. J. Math. {\bf 11}(14), 1--42 (2025)
\bibitem{St68}
Steinberg, R.: Lectures on Chevalley groups. Yale University (1968)
\bibitem{Su77}
Suslin, A.A.: On the structure of the special linear group over polynomial rings. Math. USSR-Izv. {\bf 11}, 221--238 (1977)
\bibitem{SuKo82}
Suslin, A.A., Kopeiko, V.I.: Quadratic modules and orthogonal group over polynomial rings. J. Sov. Math. {\bf 20}(6), 2665--2691 (1982)
\bibitem{Tu83}
Tulenbaev, M.S.: The Steinberg group of a polynomial ring. Math. USSR Sb. {\bf 45}(1), 139--145 (1983)

\bibitem{VaPl96}
Vavilov, N.A., Plotkin, E.: Chevalley groups over commutative rings: I. Elementary calculations. Acta Appl. Math. {\bf 45}(1), 73--113 (1996)
\bibitem{We10}
Wendt, M.: $\mathbb{A}^{1}$-homotopy of Chevalley groups. J. $\K$-Theory {\bf 5}, 245--287 (2010)
\end{thebibliography}
\end{document}